\documentclass[a4paper, 10pt, conference]{ieeeconf}      

\IEEEoverridecommandlockouts                              

\usepackage{amsmath, lipsum} 
\usepackage{amssymb}  
\usepackage{graphicx}
\usepackage{hyperref}
 
\usepackage{cite}[compress]
\usepackage{graphicx}
\usepackage{cleveref}

\usepackage{breqn}
\graphicspath{ {./figures/} }
\usepackage{float}
\usepackage{flushend} 
\usepackage{tabularx}
\usepackage{titlesec}
\usepackage{lipsum}

\titlespacing*{\section}
{0pt}{3.5ex plus 1ex minus .2ex}{2.3ex plus .2ex}
\titlespacing*{\subsection}
{0pt}{3.5ex plus 1ex minus .2ex}{2.3ex plus .2ex}
 
\title{\LARGE \bf Applications of Successive Convexification in Autonomous Vehicle Planning and Control}

\author{Ali Boyali, Simon Thompson and David Wong \\ Tier IV, Inc., Tokyo, Japan \\{ali.boyali, simon.thompson, david.wong\}@tier4.jp}
\thanks{The Authors are with Tier IV, Inc., Tokyo, Japan
        {\tt\small \{ali.boyali, simon.thompson, david.wong\}@tier4.jp}}%
 }

\begin{document}

\maketitle
\thispagestyle{empty}
\pagestyle{empty}

\begin{abstract}

In this paper, we present the application of successive convexification methods to autonomous driving problems borrowed from recent aerospace literature. We formulate two optimization problems within the successive convexification framework. Using arc-length parametrization in the vehicle kinematic model, we solve the speed planning and model predictive control problems with a range of constraints and obstacle configurations. This paper is the first systematic application of successive convexification methods from the aerospace literature to the autonomous driving problems. In addition, we show a simple application of logical state-trigger constraints in a continuous formulation of the optimization by including an evasion maneuver in the simulations section. We give details of the problem formulation and implementation, and present and discuss the results. 

\end{abstract}

\section{INTRODUCTION}
Control and planning are two key ingredients of autonomous robotic applications. Planning algorithms, in a broad sense,  provide the agents with feasible trajectories given the constraints. On the other hand, control algorithms compute and execute the instantaneous low-level control signals to follow these trajectories. Ideally, feedback controllers are desired due to their capability to account for the uncertainties inherent in the modeling and measurement processes. However, constraints relating to the state, control and other environmental factors stipulate the use of optimization methods and optimal control theories.

The motion, actuator and environmental constraint equations encountered in optimization applications are rarely linear and convex, which are the holy grail for optimization theory. Aside from some special cases, most nonlinear differential equations do not have explicit analytical solutions \cite{jordan1999nonlinear}. For optimal control problems with nonlinear dynamics, one can transform the optimization problem into nonlinear programming by using direct transcription methods \cite{betts2010practical}. Nonetheless, nonlinear differential equations are sensitive to the initial conditions along with other associated complexities which can cause the solutions to diverge \cite{mao2018successive, mao2019convexification}. 

One effective strategy to solve the nonlinear and non-convex trajectory optimization problems is to relax the problem structure by approximating the non-linear dynamics and constraints by convexification methods for real-time performance. The resulting convex optimal control problems can be efficiently solved in polynomial time thanks to the availability of state of the art convex solvers \cite{domahidi2013ecos, scs} and the increasing performance of modern computational hardware. 

A special variant of successive convexification methods (SCvx) \cite{mao2016successive, mao2018successive, mao2019convexification} have bolstered interest for solving non-convex optimal control problems in aerospace applications. There is a growing body of research showing flexible implementations of optimal control and trajectory generation for challenging motion and environment constraints. This is due to their theoretic convergence properties \cite{mao2016successive, mao2018successive, mao2019convexification} and real-time solvability by on-board computers. These algorithms can handle non-convex state and control constraints with  guaranteed convergence to nonlinear dynamics in continuous time. It has been demonstrated that sequential convexification warrants a super-linear convergence rate, owing to its constraint satisfaction property through the iterated solutions \cite{szmuksuccessive, reynolds2019dual}. The synergy of the conic optimization solvers, convexification methods and the developments in hardware afford a solution for real-time trajectory optimization with complex state and control constraints. The resultant trajectories are feasible and locally optimal.  Locally optimal and feasible trajectories are desirable due to the complications in obtaining a global optimal solution for nonlinear differential equations with real-time performance requirements. 

In this paper, we investigate the use of  SCvx algorithms \cite{mao2016successive} for autonomous vehicle motion planning and control problem formulation, implementation and simulations. This paper is the first systematic application of SCvx algorithms  from the aerospace literature to the autonomous driving problems. The SCvx algorithms are different than the classical Sequential Convex Programming (SCP) methods with additional mechanisms that ensure the strong guaranteed convergence which is not addressed in the classical SCP applications. 

In particular we formulate a full planner which gives a speed and path trajectory for autonomous driving as successive convexification problems and can be used as either a MPC controller or full planner. This improves on the frequently used iterative MPC methods by ensuring recursive feasibility and constraint satisfaction at each iteration, with proven theoretical convergence. We give the details on the use of SCVx for MPC application to path tracking and speed regulation for autonomous vehicles. Additionally, borrowing from the recent aerospace literature, we show an implementation of logical state-triggered constraints for an evasion maneuver in the simulations section. 

We structure this paper in five sections. After the introduction, we describe successive convexification briefly in Section~\ref{sec:sc}. The vehicle models and application of the SCvx methods are given in Section~\ref{sec:models}. We detail the formal definitions of the optimization problems in Section~\ref{sec:opt}. The simulation results are discussed in Section~\ref{sec:simres}. The paper is concluded in the Section~\ref{sec:conc}. 

\section{Successive Convexification - Current Theme}
\label{sec:sc}
SCvx algorithms were first presented in \cite{mao2016successive} for nonlinear motion equations and extended to handle non-convex state and control constraints in \cite{mao2017successive, mao2018successive}. The main idea is based on converting the nonlinear and non-convex equations to their convex counterparts and solving the convex sub-problems sequentially. At each iteration, the non-linear equations are linearized around the previous trajectory and non-convex constraints are transformed to convex constraints. Each of the sub-problems is solved to full optimality. Similar approaches have been proposed in the literature for autonomous vehicles. In  \cite{carvalho2013predictive, gray2013integrated}, an iterative linearization method for Model Predictive Controllers (MPC) was proposed to generate the lateral and longitudinal motion control commands for highly nonlinear vehicle dynamics equations. 

Iterative linearization and more general SCP methods work well in practice, however, the solvers might halt and give an error for infeasible solutions over the iterations even though the original problem might have a feasible solution. This behavior is called artificial infeasibility \cite{mao2016successive}. The source of this infeasibility might be coming from either linearization errors and related constraint violations or an initial infeasible trajectory. 

As an example to the infeasibility caused by linearization is the linearization of the system equations for an unrealistically too small \cite{mao2019convexification} for which no feasible input exists. 

The other common problem in convexification is the artificial unboundedness associated which is also associated with linearization. In the literature, these problem are solved by introducing a constant radius trust region for each time step on the control variables  \cite{carvalho2013predictive, gray2013integrated, alrifaee2018real}. To alleviate the constraint violation, hard constraints are converted to soft constraint on the decision variables in every constraint equations except nonlinear dynamics. However, a feasible initial solution for the system dynamics equations is still required to proceed in the iterations. Therefore, soft constraints on the boundary equations partially alleviate the infeasibility problems. 

In SCvx algorithms, these problems are handled by scheduling a trust-region and adding a virtual control signal to the equations eliminating. These simple and effective interventions prevent artificial infeasibility and artificial unboundedness \cite{mao2018successive, mao2019convexification, mao2018tutorial}. 

The trust region size is adjusted during the iterations which contributes to the guaranteed convergence property. It is also is an improvement over the constant trust region approach seen in the iterative linearization methods. Trust-region scheduling solves the artificial unboundedness. Meanwhile, the virtual control input makes the solution trajectories temporarily one-step reachable in the convex sub-problems, thus preventing the artificial infeasibility. Introducing a virtual control input to the optimization problem brings more general approach to handling infeasibility than adding a slack variable to each constraint equation.  The virtual control term, independently from the original control, acts as a soft constraint on the state equations moving the states to temporarily to a reachable region. 

On the other hand, the convergence of SCP algorithms used in the autonomous driving applications do not report a general convergence property to the nonlinear solution. One of the major contribution of the SCvx algorithms is the strong convergence property which is a relatively new development in this direction. In the SCP applications, the convergence is not usually addressed and if so, weak convergence is provided by heuristic arguments. The general strong convergence property one of the most desired aspect of solution for the real-time applications \cite{mao2016successive, mao2018successive}.  
  
One more clear distinction of the SCvx algorithms from the classical SCP methods is the integration of the system differential equations. In the SCP literature for autonomous vehicles, the equations are integrated using a constant step size, i.e using Euler integration with a single shooting propagation. In contrast, SCvx solves the integral equation between the time-steps using more accurate integration methods additionally dividing the integration into sub-intervals.  The integration step is executed in a multiple shooting manner which also contributes the general convergence property of the SCvx algorithms. The differences of the SCvx algorithms over the conventional SCP methods are summarised in the following table. 
 
\begin{table}[]
\caption{Distinctions of SCvx over General SCP Methods}
\label{tab:my-table}
\begin{tabular}{lcc}
\cline{2-3}
              & Generic SCP Algorithms & SCvx                                           \\ \hline
Trust Region  & \checkmark (constant)                     & \checkmark (scheduled -dynamic) \\ \hline
Virtual Force & -                      & \checkmark                     \\ \hline
Integration   & single step        & multi-step                               \\ \hline
Propagation   & single shooting        & multiple shooting                               \\ \hline
Convergence   & rarely addressed & strong, guaranteed                                      \\ \hline
\end{tabular}
\end{table}

The SCvx algorithm and its variants have been implemented for various aerospace problems, from powered descent \cite{szmuk2018successive, reynolds2019state} to spacecraft rendezvous \cite{malyuta2019fast}, stochastic motion planning \cite{ridderhof2019minimum, vinod2018stochastic} and many others together with the application of state-triggered constraints \cite{szmuksuccessive}. In the thesis by Szmuk \cite{szmuksuccessive}, promising computational performance is reported for a range of problem configurations. In \cite{szmuk2017convexification} agile maneuvers of a quad-rotor were demonstrated with state-triggered obstacle constraints. Integer constraints are handled with successive convexification in \cite{malyuta2019fast} for fast trajectory generation. In this paper we follow the footsteps of aerospace literature to demonstrate the application of SCvx to autonomous driving tasks. 

\section{Autonomous Vehicle Applications}
\label{sec:models}
The race for launching the first commercially viable autonomous vehicle has been an ongoing effort in the vehicle industry. Autonomous vehicles are required to handle many tasks without human supervision. Motion planning is the foremost requirement for autonomy. Real-time optimization methods similar to MPC have been prevalent for planning and control in constrained vehicle motion studies over the last decade. 

Many problems need to be addressed with suitable strategies for autonomous driving applications. For example, in highly dynamic environments, autonomous vehicles should be able to achieve safe stopping given boundary conditions in an emergency situation. If this is not possible, the planning modules should be able to find an evasion maneuver and decide another sequence of actions. These decisions should be made in a split second. Even for a simple stopping, it is obvious that, in autonomous driving, various tasks are required to be addressed. The tasks for autonomous driving include, but are not limited to, the computation of the following:

\begin{itemize}
 \item Yaw rate trajectories for yaw and roll stability during lane changing or taking a turn,
 \item Longitudinal velocity trajectories for cruise control, adjusting the car-following distance and safe stopping,
 \item Velocity and path planning for evasive maneuvers,
 \item Feasible trajectories that the low-level controllers can act on within the actuator limits
\end{itemize}

The above tasks all require fast computational modules. In this paper, we present the application of SCvx algorithms to a full plan with a velocity trajectory. We also demonstrate the use of logical state-triggers in a simple evasion maneuver. 

The successive convexification algorithms can be used for all of the above. These algorithms have been demonstrated to handle a broad range of stringent constraint requirements in real-time aerospace applications. We consider that SCvx algorithms bring improvements to autonomous driving on top of iterative MPC methods presented in \cite{carvalho2013predictive, augugliaro2012generation, carvalho2016predictive, gray2013integrated}. The improvements mainly come from the constraint satisfaction and theoretically proven guaranteed convergence. 

In the literature, various approaches to path and speed planning have been proposed with a breadth of constraint configurations and solution methods. A list and detailed comparison of various velocity planning algorithms for autonomous driving can be found in \cite{zhang2018toward}. In this study, a simple point mass kinematic model parametrized in polar coordinates is used to formulate the motion and constraint equations.  In our study, we handle a similar kind of optimization formulation with a diverse range of constraint expressions in the formulation using the kinematic vehicle models. We formulate the velocity planning problem coupled with path tracking. 

\subsection{Vehicle Model}
Under this section, we give a brief description of the vehicle models we use in the planning and low speed control problems. It is beneficial to include the vehicle models as it makes more clear the variables and constraints used in the optimization problem. Please note that, the SCvx algorithm is a general approach and admits models with any complexity.  

\begin{figure}[ht]
    \centering
    \includegraphics[width=0.3\textwidth]{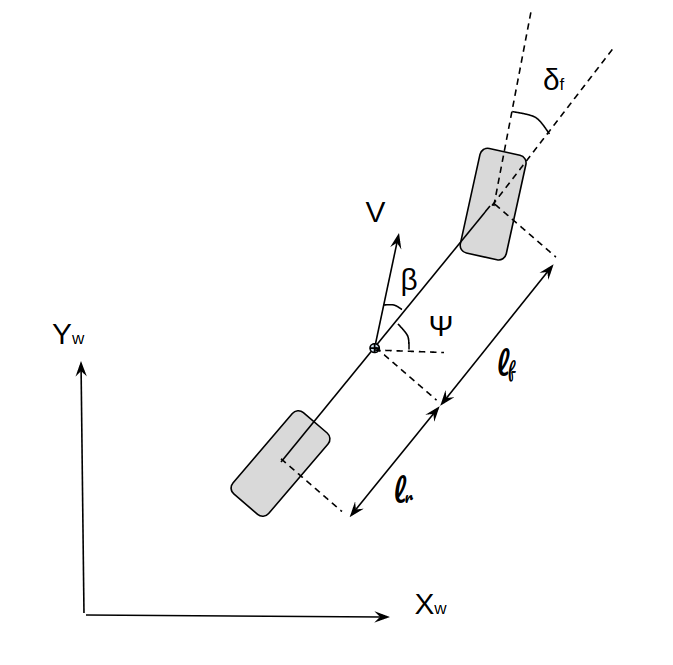}
    \caption{Single Track Vehicle Model with Side-slip Angle}
    \label{fig:single_track} 
\end{figure}

We can use the single track kinematic vehicle models with and without side-slip angle in planning simulations (see Fig. \ref{fig:single_track} and Eqs. \ref{eq:with_beta}). The first column in the equations represents a robot-car model \cite{laumond1998robot}. In this model, the rear-axle center tracks the given reference trajectory.  The second column (right) represents the same model, but considers the side-slip angle \(\beta\) \cite{kong2015kinematic, rajamani2011vehicle}. 

\begin{align} 
    \dot{X_w} &=  V\cos{(\Psi)}  & \dot{X_w} &=  V\cos{(\Psi + \beta)}  \nonumber\\ 
    \dot{Y_w} &=  V\sin{(\Psi)} & \dot{Y_w} &=  V\sin{(\Psi + \beta)} \nonumber\\
    \dot{\Psi} &= \frac{V}{l_r}\tan{(\delta_f)} & \dot{\Psi} &=  \frac{V}{l_r}\sin{(\beta)} \nonumber\\
    \dot{V} &= u_0 & \dot{V} &= u_0 \nonumber\\
    \dot{\delta_f} &= u_1 & \dot{\delta_f} &= u_1 \nonumber\\
    & &\beta &= \arctan{(\frac{l_r}{l_r + l_f}\tan{(\delta_f)})} \nonumber \\
    \label{eq:with_beta}
\end{align}

\begin{figure}[ht]
    \centering
    \includegraphics[width=0.35\textwidth]{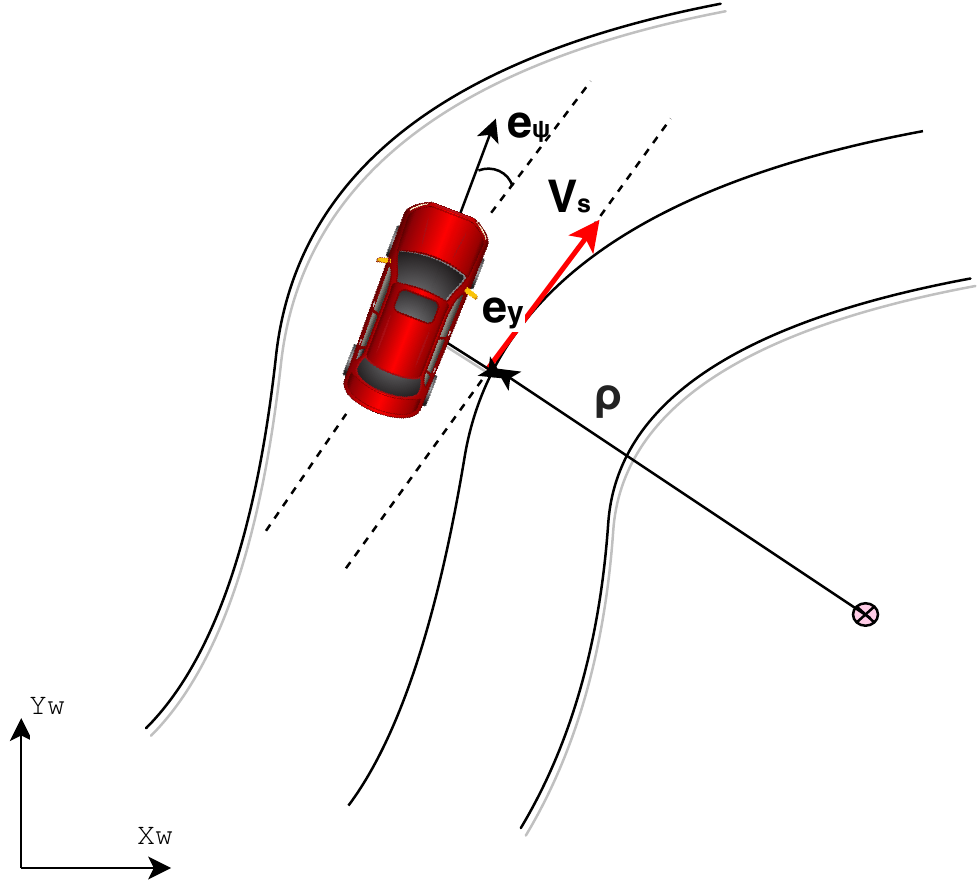}
    \caption{Road Aligned Error Coordinates}
    \label{fig:error_coord} 
\end{figure}

In Eq. (\ref{eq:with_beta}), the first three rows represent the differential equations of the global coordinates; $X_w,\; Y_w$  and heading angle $\Psi$.  The parameters $l_r,\; l_f$ mark the location of the center of gravity with respect to rear and front axles. The acceleration $u_0 = \dot{V}$ and the steering rate $u_1 = \dot{\delta_f}$ are the control inputs to the models. To capture the affect of tire deflection on the dynamics, the side-slip angle $\beta$ is included in the models which can be used for path tracking at lower speeds. 
We transform the model into error dynamics by expressing the equations in the error coordinates  as in \cite{gray2013integrated, carvalho2016predictive, qian2016model} (Fig. \ref{fig:error_coord})).

In addition, we use arc-length parametrization instead of time in the implementations. The travelled distance $s$ is used as the integrating variable instead of time. The expression of this transformation for the side-slip angle model is given in Eqs. \eqref{eq:xprime} and \eqref{eq:arc-length}. A similar transformation can be obtained from the left column (robot-car) equations of Eq. \ref{eq:with_beta}. By these modifications, not only can the obstacle constraints along the path be expressed linearly, but also the trajectory problem is converted into a path planning problem. The arc length parametrization brings another advantage in terms of the planning horizon which can be chosen arbitrarily long in the SCvx algorithms. This makes the problem a fixed time optimal control problem without introducing an additional non-linearity due to the time scaling. 

\begin{align}
    \dot{x(t)} &= f(x(t), u(t)), \nonumber & x(s)^{\prime} &= \frac{dx}{dt}\frac{dt}{ds} = f(x(t), u(t)) \frac{1}{\dot{s}} \nonumber\\
    \label{eq:xprime}
\end{align}

After the transformations and the arc-length parametrization, the final set of the equations in the error coordinates takes the form of 

\begin{equation}
    x^{\prime}(s) = F(x(s), u(s), \kappa(s))
    \label{eq:arc-length}
\end{equation}
 
\noindent with the states:
\begin{equation}
\begin{aligned}[c]
\dot{s} &= \frac{1}{1 -\kappa e_y} V \cos{(e_\Psi + \beta)}\\ 
e_y^{\prime} &= V \sin{(e_\Psi + \beta)} \frac{1}{\dot{s}} \\
 V^{\prime} &= u_0 \frac{1}{\dot{s}} \nonumber \\
\end{aligned}
\;
\begin{aligned}[c]
\Psi^{\prime} &=  \frac{V}{l_r}\sin{(\beta)} \frac{1}{\dot{s}} \\
 e_\Psi^{\prime} &= \Psi^{\prime} - \kappa(s)\\
 \delta_f^{\prime} &= u_1 \frac{1}{\dot{s}}  
\end{aligned}
\end{equation}

The time $t_s^{\prime}= \frac{1}{\dot{s}}$ can also be included in the states to recover the time dependent trajectories. In the equations, the road curvature $\kappa(s)$ enters the model as the reference parameter. It can be taken from off-line look-up tables or approximated from sensor measurements. 

\subsection{Convex Sub-problems}
The successive convexification method is built upon the convex sub-problems which arise from linearization around the previously computed or predicted trajectories. We can express the general structure of the convex sub-problems as in the following form;

\begin{align}
    \min_{u(s)}  \quad & J(s(t_0), s(t_f), x(s), u(s), \kappa(s))  \\ 
    \textrm{s.t.} \quad  & h(x(s), u(s), \kappa(s)) =0  \label{eq:equality_constraint}  \\
     & g(x(s), u(s), \kappa(s)) \leq 0 \label{eq:inequality_constraint} 
\end{align}

\noindent where the equality and inequality constraints are given by \eqref{eq:equality_constraint} and \eqref{eq:inequality_constraint} and the boundary conditions of the integration are defined as $s \in [s(t=0),\; s(t=t_f)]$. The final distance is decided depending on the problem configuration.  

\subsection{Linearization, Discretization and Scaling}

The linear time-varying model of the vehicle motion is approximated by the first order Taylor expansion as;
\begin{align} 
    x^{\prime}(s + ds) & \approx A(s) x(s)+B(s) u(s)+F(s) + w(s)  \nonumber\\ 
    A(s) &:=\left.\frac{\partial F}{\partial x}\right|_{\bar{z}(s)}  \nonumber\\ 
    B(s) &:=\left.\frac{\partial F}{\partial u}\right|_{\bar{z}(s)}  \nonumber\\ 
    F(s) &:= F(\bar{x}(s), \bar{u}(s),\bar{\kappa}(s))  \nonumber\\
    w(s) &:=-A(s) \bar{x}(s)-B(s) \bar{u}(s)  \nonumber \\
    \label{eq:discretization}
\end{align}

\noindent where the bar notation represents the computed trajectory coordinates obtained by the previous solution ($\bar{z(s)}$) at which the linearization takes place with: 

$$\bar{z}(s):=\left[\bar{s}\;, \bar{x}^{T}(s)\;, \bar{u}^{T}(s)\right]^{T} \text { for all } s \in [0, 1].$$ 

The real arc-length $s_{real}  \in [s_{t=0}, s_{t=t_f}]$ is re-scaled into the interval $s \in [0,\;1]$. In the speed planning study, the final point of the arc-length and the number of intervals are fixed. 

The equality and inequality constraints are convexified if they are not convex in a similar manner: 

\begin{equation}
    h_{i}(z(s)) \approx h_{i}(\bar{z}(s))+\left.\frac{\partial h_{i}}{\partial z}\right|_{\bar{z}(s)} \delta z(s)
\end{equation}

\begin{equation}
    g_{i}(z(s)) \approx g_{i}(\bar{z}(s))+\left.\frac{\partial g_{i}}{\partial z}\right|_{\bar{z}(s)} \delta z(s)
\end{equation}

\noindent with the definition of $\delta z(s):=z(s)-\bar{z}(s)$. 

The performance of various discretization methods for successive convexification was reported in \cite{szmuksuccessive, malyuta2019discretization}. According to this study, the First Order Hold (FOH) has the lowest computational time among others (zero order hold, classical Runge-Kutta and pseudo-spectral methods). Therefore, we used the FOH method to discretisize the state trajectories and the control using the fundamental matrix solution to the ordinary differential equations. In the FOH discretization, the parameters and input the models are interpolated between the start and end points of each integration interval. In this case, the discrete equations take the following form \cite{szmuksuccessive, malyuta2019discretization} with the intervals $s_k = \frac{k}{K-1}$ and $k \in[0, K]$ for the each step $k$: 

\begin{dmath}
    x^{\prime}(s) =A(s) x(s)+\lambda_{k}^{-}(s) B(s) u_{k}+\lambda_{k}^{+}(s) B(s) u_{k+1}  + F(s) + w(s) 
\end{dmath}

$\forall s \in\left[s_{k}, s_{k+1}\right]$ and the interpolating coefficients:

\begin{align}
    u(s) &=\lambda_{k}^{-}(s) u_{k}+\lambda_{k}^{+}(s) u_{k+1}, & \forall \tau \in\left[s_{k}, s_{k+1}\right] \nonumber\\
    \kappa(s) &=\lambda_{k}^{-}(s) \kappa_{k}+\lambda_{k}^{+}(s) \kappa_{k+1}, &  \nonumber\\
    \lambda_{k}^{-}(s)& :=\frac{s_{k+1}-s}{s_{k+1}-s_{k}} &  \lambda_{k}^{+}(s):=\frac{s-s_{k}}{s_{k+1}-s_{k}} \nonumber\\
\end{align}

\noindent with the fundamental matrix solution \cite{antsaklis2006linear, hespanha2018linear} to ODEs:

 \begin{equation}
     \Phi_{A}\left(\xi, s_{k}\right)=I_{n_{x} \times n_{x}}+\int_{s_{k}}^{\xi} A(\zeta) \Phi_{A}\left(\zeta, s_{k}\right) d \zeta
 \end{equation}
 
 \noindent where $n_x$ is the number of states in the equations. Using the properties of the fundamental matrix (Theorem II.1 in \cite{malyuta2019discretization}) one can arrive at the non-homogeneous solution  of the Linear Time Varying (LTV) system equations as:
 
  \begin{align}
    x_{k+1} &=A_{k} x_{k}+B_{k}^{-} u_{k} + B_{k}^{+} u_{k+1} + F_{k} + w_{k} \\
    A_{k} &:=\Phi_{A}\left(s_{k+1}, s_{k}\right) \\ 
    B_{k}^{-} &:=A_{k} \int_{s_{k}}^{s_{k+1}} \Phi_{A}^{-1}\left(\xi, s_{k}\right) B(\xi) \lambda_{k}^{-}(\xi) d \xi \\
    B_{k}^{+} &:=A_{k} \int_{s_{k}}^{s_{k+1}} \Phi_{A}^{-1}\left(\xi, s_{k}\right) B(\xi) \lambda_{k}^{+}(\xi) d \xi \\ 
    F_{k} &:=A_{k} \int_{s_{k}}^{s_{k+1}} \Phi_{A}^{-1}\left(\xi, s_{k}\right) F(\xi) d \xi \\
    w_{k} &:=A_{k} \int_{s_{k}}^{s_{k+1}} \Phi_{A}^{-1}\left(\xi, s_{k}\right) w(\xi) d \xi  
 \end{align} 

These matrices are used as the inputs to the solvers to compute the motion equation constraints in the implementation. 

One more subtlety in the numerical optimization is scaling. In general, the states in the equations do not have a similar range of magnitude. This may introduce inconsistencies into the numerical solution. To prevent such problems during the optimization, it is a common practice to normalize the states. One way to do scaling is with a linear transformation. We applied the following transformation to all of the  states $x$ and the inputs $u$ \cite{reynolds2019dual, gill1981practical, ross2018scaling};  

\begin{gather}
    x =D_{\hat{x}} \hat{x}  + C_{\hat{x}} \label{eq:scalingx}\\
    u =D_{\hat{u}} \hat{u}  + C_{\hat{u}} \label{eq:scalingu}\\
    \nonumber
\end{gather}

\noindent where $\hat{x}$ and $\hat{u}$ are the normalized state and control variables. The solver seeks a solution in the normalized variables denoted by the hat notations in \eqref{eq:scalingx}, \eqref{eq:scalingu}. The scaling coefficient matrices $D_{\hat{x}},\; D_{\hat{u}} $ and the centering vectors $C_{\hat{x}}, \; C_{\hat{u}}$ can be computed from the maximum and minimum boundary values of the variables. 

\subsection{Virtual Force for Artificial Infeasibility and Trust Region for Artificial Unboundedness}
A virtual force vector $\nu_k \in \mathcal{R}^{n_x}$ is added to the system equations as an input to prevent artificial infeasibility, making the states one-step reachable: 
\begin{equation}
    x_{k+1}=A_{k} x_{k}+B_{k}^{-} u_{k}+B_{k}^{+} u_{k+1}+w_{k}+\nu_{k}    
\end{equation}

The $\ell_1$ norm is used in the cost function of the virtual force to promote sparsity \eqref{eq:l1}. A high penalty weight $w_{\nu}$  is assigned in the cost for the virtual force so that the solver uses it only when necessary.  The following cost is added to the optimization objective cost function:

\begin{equation}
    J_{v}(\nu):=w_{\nu} \sum_{k \in \overline{\mathcal{K}}}\left\|\nu_{k}\right\|_{1}
    \label{eq:l1}
\end{equation}

As the linearization step might introduce unboundedness, the search space in the optimization variables $[x-x_k,\;u-u_k]$ are bounded by a trust region. Two forms of trust region formulation are practised in the literature. One can add a hard trust region constraint into the optimization formulation or a trust region cost can be added to the optimization cost using the soft trust region method \cite{benedikter2019convex, szmuksuccessive}. We used the former one: 

\begin{equation}
     \left\lVert \delta x_{k}  \right\lVert_1 + \left\lVert \delta u_{k}  \right\lVert_1 < \rho_{tr} 
\end{equation}

\noindent where $\delta x_{k} := x - x_k,\;\delta u_{k} := u - u_k,$ and $\rho_{tr}$ is the trust region radius which is scheduled depending on the accuracy of the linear approximation in the model. The details of the trust region scheduling algorithm are given in \cite{mao2018successive, mao2018tutorial}. 

\subsection{State-Triggered Constraints}
In some applications, we may need some constraint when a certain criterion is met in a continuous optimization problem. An evasion manuever can be given as an example to this case. If a collision is obvious on the way through a trajectory plan we can reformulate the alternatives in a continous optimization to make the autonomous agent to perform another maneuver. In general, embedding logical constraints require stringent treatment and and mixed integer programming is involved. There are specialized optimization software for this purpose. Szmuk et al. in \cite{szmuk2018successive_trig} proposed state-triggered constraint formulation in their successive convexification applications and show its use in different optimization problems improving the method in \cite{reynolds2019dual, reynolds2019state, szmuk2017convexification}. The underlying structure of continuous state-triggered constraint formulation is very similar to linear complementarity problem \cite{cottle1992linear}.  Assume that we want to enforce some constraints equations $c(z) = 0$ given some conditions $g(z) \le 0$ with a variable $z \in \mathcal{R}^{n_z}$ belongs to parent optimization problem.  The formal definition follows as:

\begin{align}
    \sigma & \geq 0  \nonumber\\  
    g(z)+\sigma & \geq 0  \nonumber\\ \nonumber
    \sigma c(z) & \leq \label{eq:trig_c} 
\end{align}

\noindent where a new slack variable $\sigma \in \mathcal{R_{++}}$ is introduced. The equations \ref{eq:trig_c} admits an analytical solution with  $\sigma^{*}:=-\min (g(z), 0)$ (please see \cite{szmuk2018successive_trig, reynolds2019state,reynolds2019dual } for the details). Thus the state-triggered constraints in a continuous optimization can be formulated as $h(z):=-\min (g(z), 0) \cdot c(z) \leq 0$. The variable $z$ can take any parameter from the previously solved trajectories in the successive convexification. In our demonstration, we devised a simple evasion maneuver in which if an autonomous vehicle cannot reach a desired speed (let's say 1 m/s) at terminal point of the trajectory, we request the vehicle to avoid from a virtual obstacle at the terminal location. 

\section{Problem Formulations in Planning and Control}
\label{sec:opt}
We applied the successive convexification methods to a full planning (speed and path) and MPC path tracking problems. However, for the sake of clarity, we narrow down the application to only planning problems in the paper. The method can be extended to MPC formulation by including a receding horizon. Expressing the path region and obstacle constraints as well as the terminal conditions are straightforward with arc-length parameterization. Therefore, we used arc-length parametrization in the solution. 

\subsection{Application of SCvx for Speed Planning}
We consider the scenarios of reducing or increasing the vehicle speed within a given travelling distance while the car is following a curved path, with state and actuator constraints.  The main objective is to arrive at a speed value at the end boundaries of the solution. For example, if there is a pedestrian in a certain distance, the autonomous vehicle must be able to stop before the pedestrian location.  One can make the problem more complex by adding an obstacle avoidance constraint in the speed planning problem as well. Thus, the problem becomes full planning problem in which the path avoiding the obstacle needs also be generated along with the speed plan. 

We simulated various scenarios using the vehicle models given by Eq. \eqref{eq:with_beta}. In planning algorithms, less complex models are preferred in terms of computational speed. Therefore, we give the results of the planning for only the vehicle model without the side-slip angle formulation given in the first column of Eq. \eqref{eq:with_beta}. The arc-length parametrization introduces singularity around the point $\dot{s} = 0$ in the equations as the final speed is required to be zero for the stopping problems. For planning applications this singularity can be ignored, however when recovering the time dependent control signal it can be avoided by defining a lower bound of speed other than zero such as $V_{final} = 0.5 \; [m/s]$. From this speed value, we assume that the vehicle can be stopped safely. 
 
The formal definition of the speed planning problem given the boundary conditions is formulated in the discrete form with $s := s_{[1:K]}$ as follows: 

\begin{align}
    \min_{u(s_{1:K})} \quad & J(s(t_0), s(t_f), x(s), u(s), \kappa(s))  \\ 
    \textrm{s.t.} \quad  & h_i(x(s), u(s), \kappa(s)) =0   \\
     & g_i(x(s), u(s), \kappa(s)) \leq 0  
\end{align}
where the controls are acceleration and steering rate and the cost functionals are given by; 

\begin{align}
    J(s, x, u, \kappa) & = J_{e_y} + J_{e_\Psi} + J_{\nu} +  J_{jerk}  + J_{u} + J_{N = s_f}   \nonumber\\
    J_{e_{y, \Psi}} &= w_{e_{y, \Psi}} \left\lVert e_{_{y, \Psi}{[1:K]}} \right\lVert_2  & \nonumber \\
    J_{jerk} &= w_{jerk} \left\lVert \Delta u_{0_{[1:K-1]}} \right\lVert_2  \nonumber \\
    J_u &= w_{u_0} \left\lVert u_{0_{[1:K]}} \right\lVert_2 + w_{u_1} \left\lVert u_{1_{[1:K]}} \right\lVert_2 \nonumber\\
    J_{N=s_f} &= w_{N} \left\lVert V_{[K]} - V_{N=s_f} \right\lVert_2  \nonumber \\
    \label{eq:discrete_opt}
\end{align}

\noindent where $w_{(.)}$ represents the weights of the associated variables. In the implementation, all the weights are set as unity except the virtual force weight which must be set sufficiently high. We included longitudinal jerk cost to reduce the rapid changes in the acceleration considering the comfort requirements. The jerk cost is formulated as the two-norm of the difference of successive acceleration inputs. Since the dynamics are parametrized by the arc-length, this jerk expression is called  pseudo-jerk for the longitudinal motion. The last term $J_{N=s_f}$ captures the terminal constraint for the desired speed at the end of the planning horizon. The desired speed can also be enforced in the constraint equations. Although the the range of controls are constrained in the optimization, one can define additional cost to further reduce the amount of control effort exerted into system. The lateral deviation and deviation of heading angle $J_{e_{y, \Psi}}$ can be enforced to have the vehicle to follow the prescribed path while performing the main objectives (i.e while trying to stop at the given horizon or avoiding an obstacle). 

We formulated the cost terms using the norms of the variables as recommended in the CVX optimization library \cite{CVXeliminating} as we use the Second Order Conic Programming (SOCP) solver ECOS \cite{domahidi2013ecos} for solving the problems. 

We enforced the following state and control constraints;

\begin{equation}
\begin{aligned}[c]
 V   & \leq V_{max}, \; V \in \mathcal{R}_{++} \\
 \lvert \delta \rvert & \leq 27 \; [deg] &  \\
 x(s_0) & = x_{initial} \\
 x(s_f) & = x_{final} \\
\end{aligned}
\;
\begin{aligned}[c]
\lvert \dot{\delta_f} \rvert & \leq 60 \; [deg/sec] &  \\
    \left\lVert
    \begin{bmatrix}
        a_{x}, & a_{y}\\
    \end{bmatrix} \right\rVert_2 & \leq \mu g  & \label{eq:friction0}\\ 
    u(s_0) &= u_{initial} \\
    u(s_f) &= u_{final} \\
    e_y(s) & \leq p(s) \\ 
\end{aligned}
\end{equation}

In the constraint equations, $e_y(s) <= p(s)$ can be used for the boundary of the road or to define the obstacle avoidance distance along the path. It can be a piece-wise linear function or a nonlinear function that can be convexified. In the simulations, we locate a virtual obstacle along the center line the vehicle tries to follow. When the vehicle reaches the obstacle regions $20:24^{th}$ discretization region (corresponding the region between 25 and 30 meters) it put a distance to this region minimum 0.5 meters from the right.     

The acceleration constraint, $\rVert [a_{x}, \;  a_{y}] \lVert_{2}\leq \mu g$, captures the friction circle constraint (Kamm’s circle) with the tire road adhesion coefficient $\mu \in [0, \; 1]$. Once again, this is a simplified assumption for capturing the tire force limit and used as demonstration only. The aim is not to formulate a full fledged vehicle model with tires but show the use of SCvx algorithms in the autonomous driving domain. One can design models and problems serving the objective in their design.  

The rate constraints in time and arc-length parametrization can be converted in either direction by simple algebraic manipulations when forming the constraint equations by affine relations. 

Given the formal description of the problem, speed profile trajectory optimization can be described as having the vehicle reach a terminal speed at the end of a given road section, while following a prescribed path and respecting the state and control constraints.

One important point in the velocity planning is to guess an initial solution. The initial solution can be prepared by generating an input sequence and integrating the equations for the given inputs, or interpolating the state and control variables between the boundary conditions. Due to the virtual force in the model, the predicted trajectories are not required to be feasible. In our simulations, in addition to linear interpolation, we initialize the steering angle related variables from the Ackermann steering geometry as the road section and curvature are all known. 


\subsection{Use of State-Triggered Constraints for an Evasion Maneuver}
The continuous formulation of state-triggered constraints allows us to formulate logical constraints in a continuous optimization problems. Otherwise, one would use specialized software suit to be able to solve logical constrained problems. Other advantage of the state-triggered constraints is the reduction of the constraint variables in the problem. For example, if there will be some region in a state-space domain where the constraints are not necessary to enforce can be defined by these framework. There is no need to define a constraint equation along the path where the constraints are inactive. 

We show the use of state-triggers with a simple evasion maneuver. Similar to full planner example given above, we set a terminal distance 50 meters, and initial speed in such a way that the vehicle cannot reach the final desired speed from the given speed. For this example, we set the initial speed to 25 [m/s] and the final speed 0.5. From the initial speed defined, the car cannot reach the final speed exactly and its speed is the final distance is around 2 [m/s]. The only difference from the obstacle free scenario presented in the previous section is the initial speed. We increased the initial speed from 20 to 25 [m/s] so that the car cannot stop in the given range. The state trigger constraint $g(z) = [V_{final} \ge 1]$ is defined which states that if the vehicle's terminal speed is greater then a one meter/second. If this condition is encountered during the solution, we requested the vehicle to start an evasion maneuver 2.5 meters before the final point (last two discrete points in the simulations).  The constraint conditions therefore becomes $c[z]= [e_{y{N-2:N}}  \ge 1]$ putting a distance a minimum of one meters from the left in the last two discrete points of the solution trajectory.  
 
\section{Simulation Results} 
\label{sec:simres}
 
We solved the speed planning problem with several constraint configurations. We enforced the desired terminal speed in the objective function. The final arc-length can be chosen arbitrarily (depending on the problem configuration). As an example, if there is a pedestrian 50 meters ahead and the car is requested to stop within this range before the pedestrian, the total arc-length in this scenario can be chosen as 50 meters. 
We observed convergence up to 100 meters total arc-length with 25-50 discretization steps. We present the result of two scenarios. The initial speed and the road-tire friction coefficient $\mu$ are set to 20 [m/s] and 0.6 respectively for both of the scenarios. In the first, the vehicle is requested to slow down to 0.5 [m/s] within 50 meters. In the second scenario, the vehicle is requested to slow down and reach at a terminal speed similar to the first one but with an additional obstacle constraint. A virtual obstacle is located between the 20-24 steps of the discretization out of 40. In this case, the vehicle must satisfy two objectives; slowing down to a terminal speed and avoiding an obstacle while doing so.  We formulated obstacle avoidance to ensure the vehicle pass the obstacle region from left or right and put an distance to the obstacle at least 0.5 meters. Please note that, this is a simple path and velocity plan demonstration. In the planning algorithms using simple models is a custom assuming that the plan can be refined and passed to more accurate controllers and planners. Therefore, for demonstration the algorithms we simply ignore the vehicle geometry. 

\begin{figure}[h]
    \centering
    \includegraphics[width=0.5\textwidth]{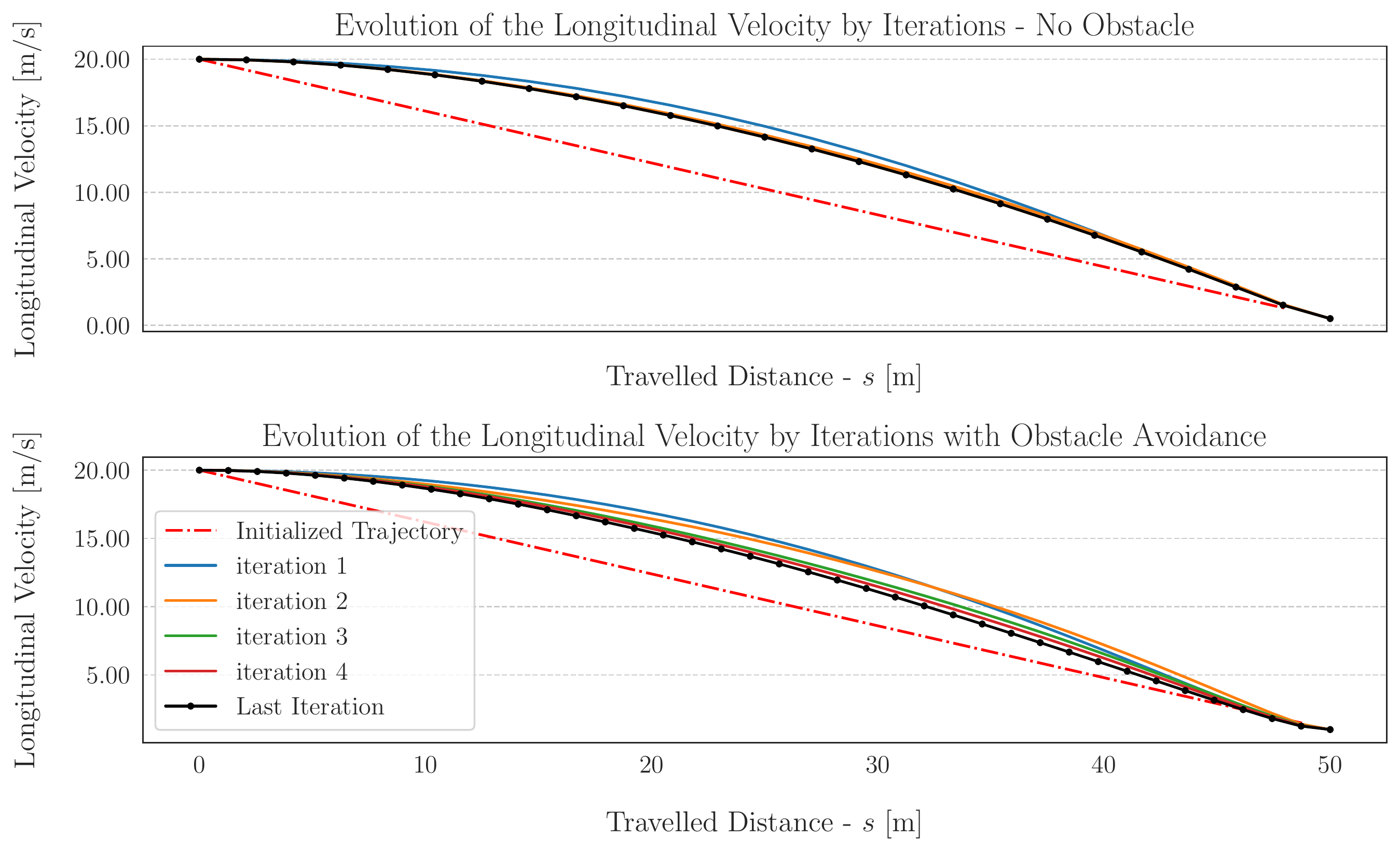}
    \caption{Evolution of Speed Trajectories with and without Obstacle Avoidance}
    \label{fig:vx_vp} 
\end{figure}

\begin{figure}[h]
    \centering
    \includegraphics[width=0.5\textwidth]{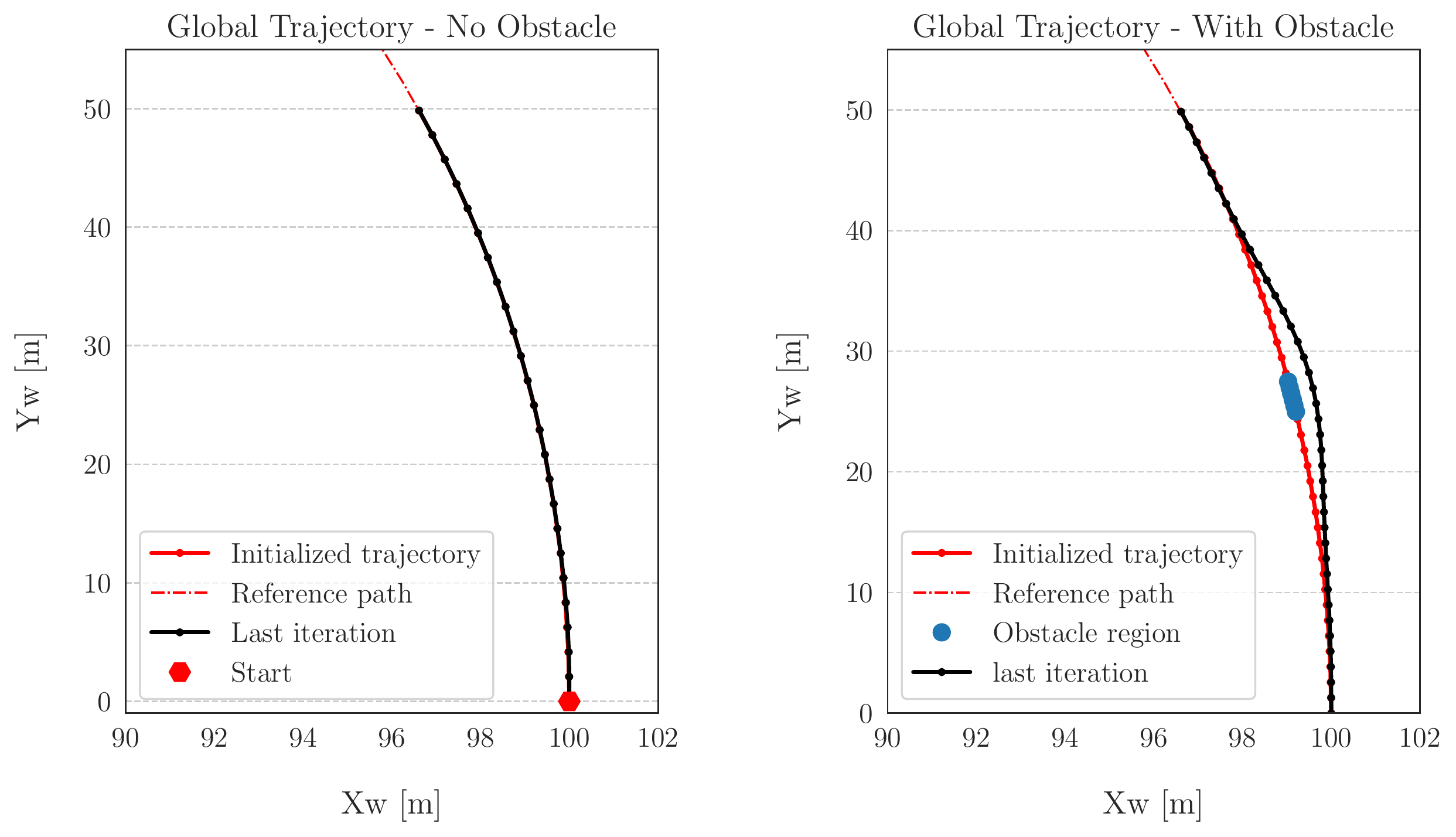}
    \caption{Evolution of Paths with and without Obstacle Avoidance}
    \label{fig:path_vp} 
\end{figure}

Fig. (\ref{fig:vx_vp}) shows the evolution of the vehicle longitudinal speed on the road section for the two scenarios while the corresponding paths are given in Fig.(\ref{fig:path_vp}) 

In Fig.(\ref{fig:vx_vp}) the initial trajectories are plotted in red with a dashed line, whereas the final trajectories with a dot-marked continuous line in black. As shown in the figure, the trajectories converge a solution. It takes more iterations to converge in the obstacle case. 


The resolution of the trajectory, the number of discretization steps and the final distance are problem dependent. We tested different final distance and number of discretization steps and were able to solve the planning problem with several combinations. One can use this flexibility to further improve the solutions by first solving with a low resolution and improving the optimization mesh if necessary. 

\begin{figure}[h]
    \centering
    \includegraphics[width=0.48\textwidth]{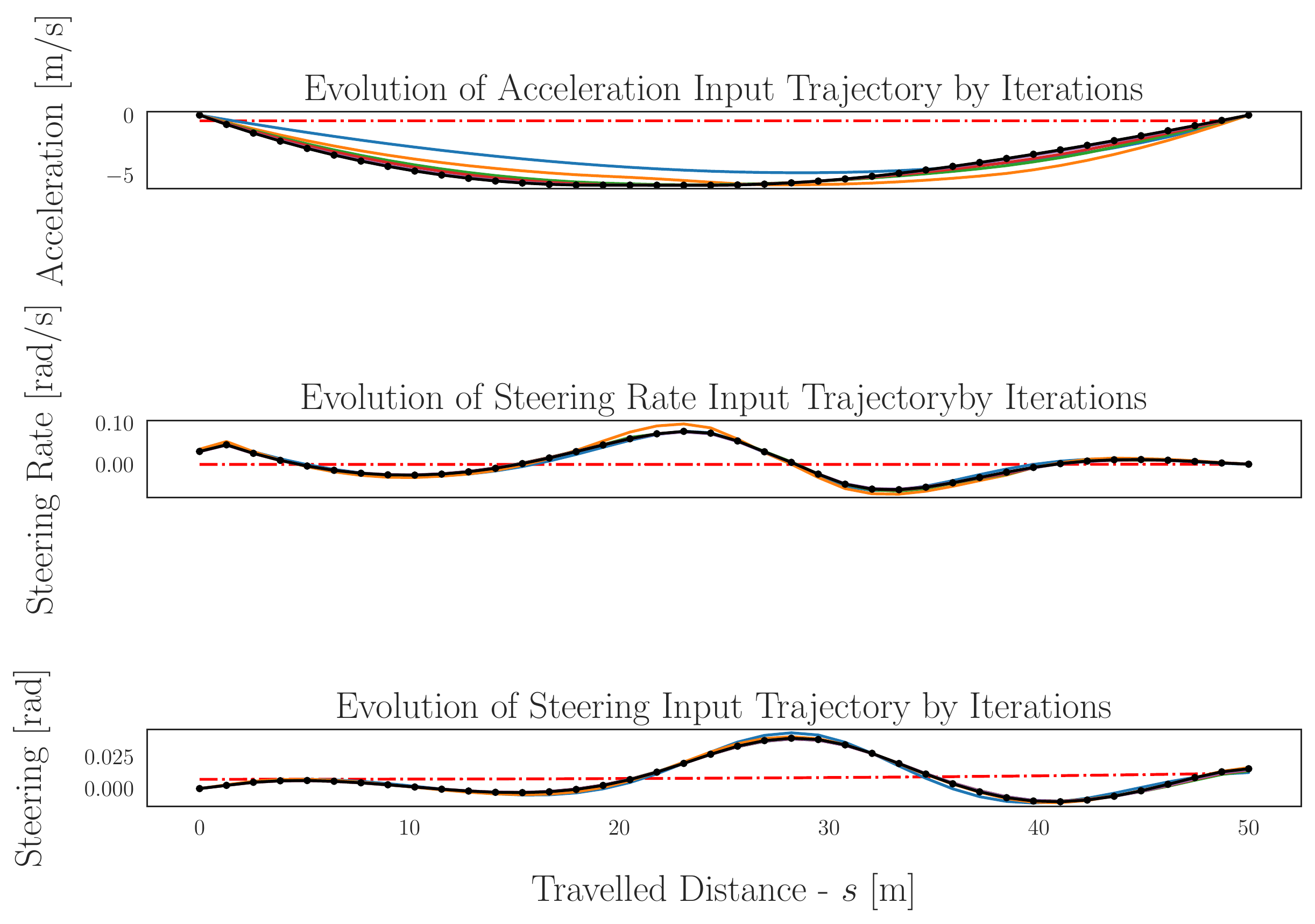}
    \caption{Evolution of Controls for Speed Planning Problem}
    \label{fig:inputs_vp} 
\end{figure}

We give the computed control inputs for the velocity planning with obstacle problem only in Fig. (\ref{fig:inputs_vp}) due to limited space. The first rows of the figure presents the evolution of computed acceleration inputs to the models. In the middle row, the steering state variable is shown. The last row is the steering rate control input to the model.  All the input constraint equations are satisfied as expected. 

\begin{figure}[H]
    \centering
    \includegraphics[width=0.4\textwidth]{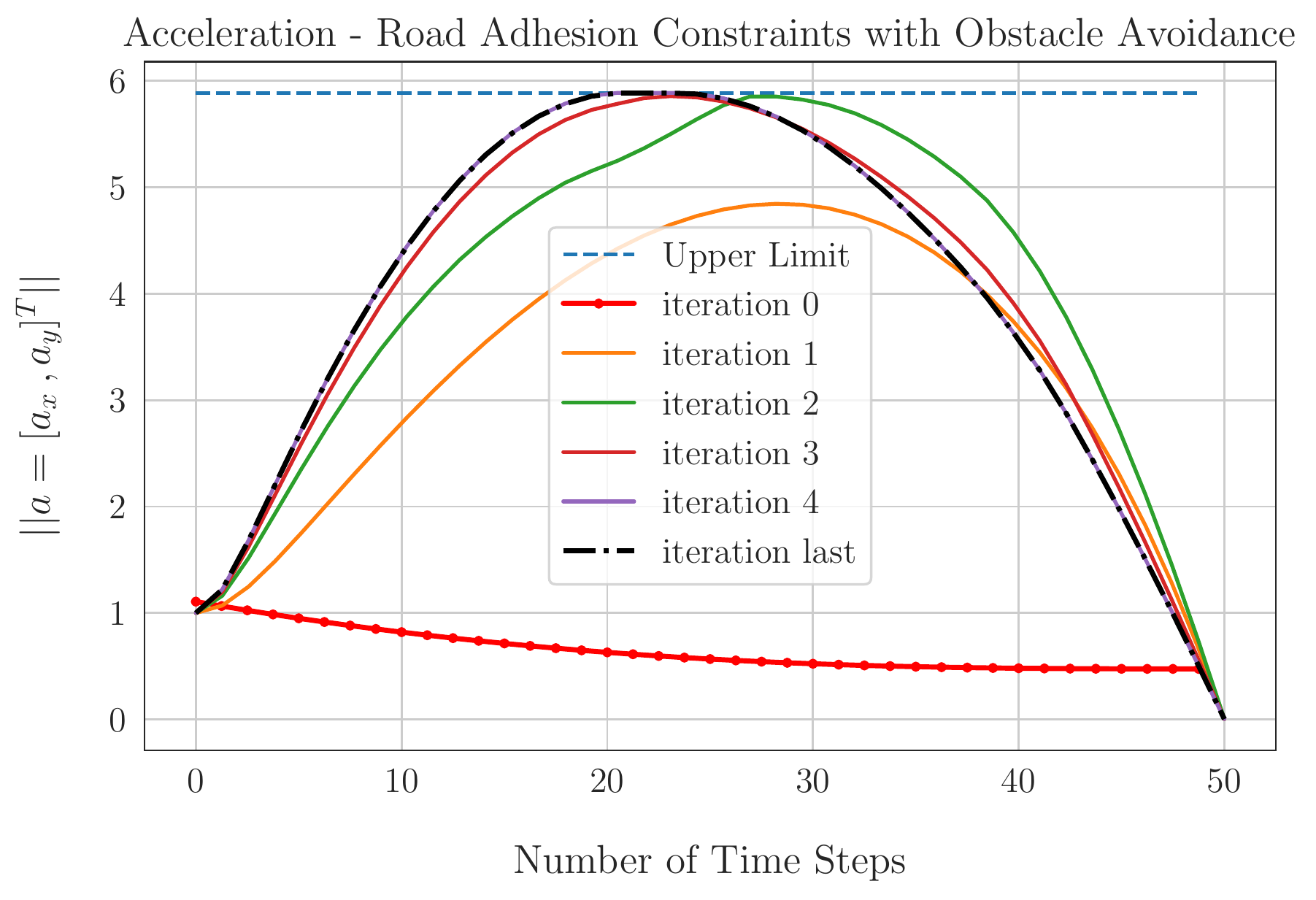}
    \caption{Friction Circle Constraint - with Obstacle}
    \label{fig:acc_vp} 
\end{figure}

Constraint satisfaction is more clear in in Fig. (\ref{fig:acc_vp}) which shows the the magnitude of the combined acceleration vectors. The solution here respects the friction circle constraint which is the dashed-line in the figure.  
 
\begin{figure}[H]
    \centering
    \includegraphics[width=0.4\textwidth]{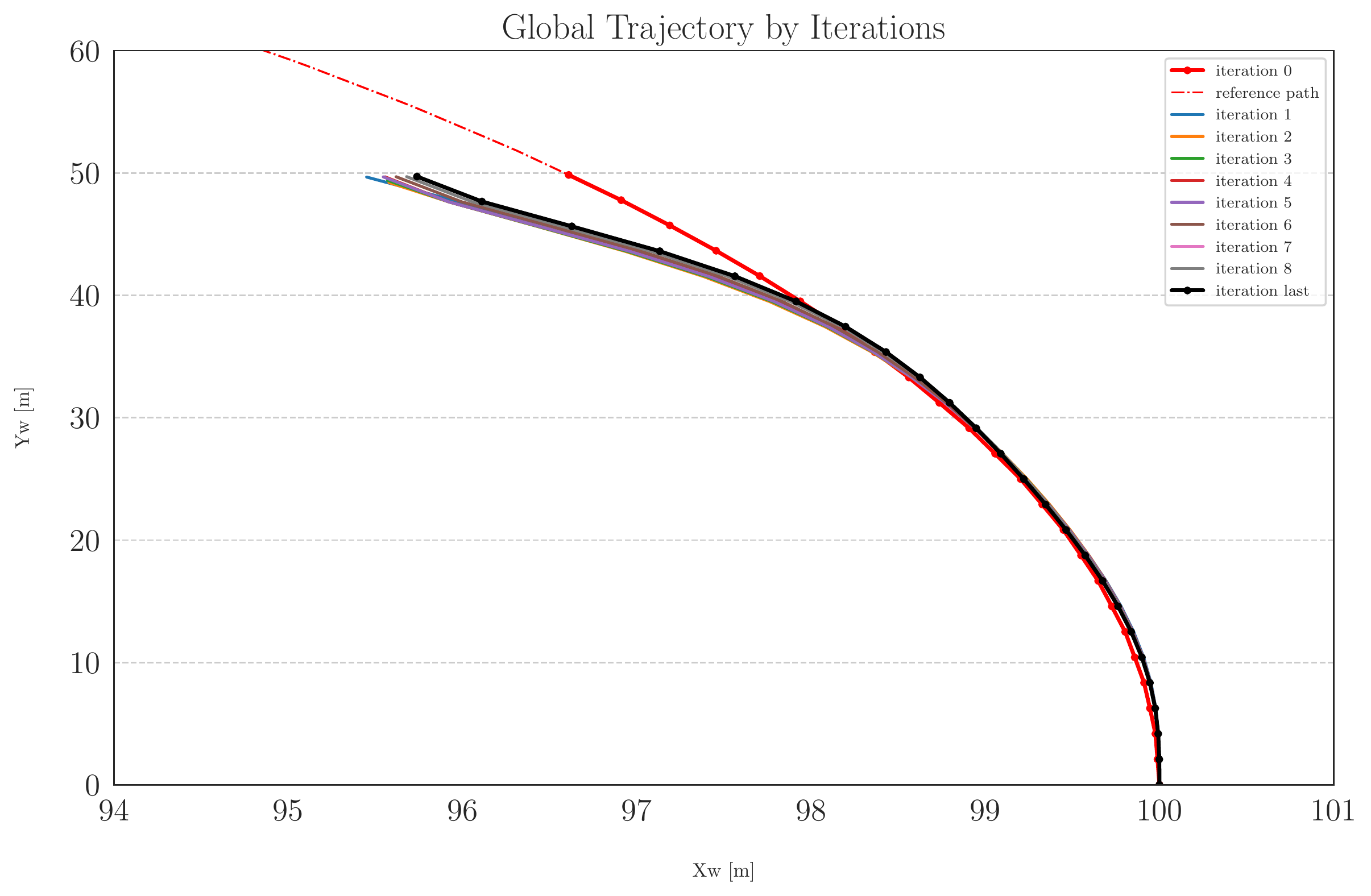}
    \caption{Global Trajectory - Result of Active State-Triggered Constraint}
    \label{fig:evasion} 
\end{figure}

The simulation results for the state triggers is shown in the Fig. (\ref{fig:evasion}). 

In the figure, it is shown that, the vehicle put a distance with a minimum of one meters from left to a virtual obstacle located at the end of the path. Since the solver cannot find a speed trajectory reaching the desired speed (0.5 m/s) at the terminal location. In this case, the state-triggered constraint is activated and the algorithm gives a solution to the second optimization problem to the evasion maneuver by which the car avoids a virtual obstacle located at the terminal point passing the obstacle from left by putting a specified distance.

\section{DISCUSSION and CONCLUSIONS}
\label{sec:conc}
In this paper, we presented the application of SCVx methods to vehicle models for velocity planning and evasion maneuver design with a state-trigger. The vehicle models are arc-length parametrized to flatten the coordinates of the curved paths. In the planning algorithm, simplified models i.e point-mass kinematic model are used to find a feasible trajectory. In the spirit of this custom, we haven't used a complex model to formulate the planning problems but we kept some complexity at a certain level to formulate realistic situations. 

We used uniform discretization intervals and demonstrated the application of SCvx algorithms to the vehicle motion planning and control problems. The SCvx methods are algorithmic methods that bring about the possibility to solve optimal trajectory generation problems in real-time. The SCVx methods bring improvements in convergence by constraint satisfaction on top of the iterative MPC algorithms in the literature. There have been a wide range of methods that facilitate convergence and increase accuracy. Among them are the well-known pseudo-spectral methods which integrate the trajectories with non-uniform discretization points. One can further improve the application of SCvx making use of the pseudo-spectral method if accuracy and precision is of primary interest. 

We implemented the problem structure in the Python environment as it is easier to manipulate the codes, to visualize the results than implementing the simulation in C++. The computational performance of the algorithms in Python is not suitable for real-time implementation for our simulations. The computations take  2.5 - 3 secs for the planning algorithms. It might be a sound computational time, however our initial C++ implementation for path tracking provides a solution at 40Hz. We did not make use of strategies such as  customized solvers \cite{dueri2014automated}, or parallelizing the discretization steps to improve computation times in this study. 

The SCvx algorithms together with the listed computational strategies are promising and viable solution for autonomous vehicle planning and control applications. In the future, we plan to experiment with these techniques on a real car and develop more flexible obstacle avoidance constraint implementations. We will release the implementation of the algorithms under Autoware open source software repository \cite{kato2018autoware, aw:autoware_repo}. 




\section*{ACKNOWLEDGMENT}

We would like to thank to Sven Niederberger for sharing his implementations of SCvx algorithms in both Python and C++ on github \cite{steven}. We are inspired from his code templates. 

\bibliographystyle{ieeetr}
\bibliography{bibliography}

\begin{thebibliography}{10}

\bibitem{jordan1999nonlinear}
D.~W. Jordan and P.~Smith, {\em Nonlinear ordinary differential equations: an
  introduction to dynamical systems}, vol.~2.
\newblock Oxford University Press, USA, 1999.

\bibitem{betts2010practical}
J.~T. Betts, {\em Practical methods for optimal control and estimation using
  nonlinear programming}, vol.~19.
\newblock Siam, 2010.

\bibitem{mao2018successive}
Y.~Mao, M.~Szmuk, X.~Xu, and B.~Acikmese, ``Successive convexification: A
  superlinearly convergent algorithm for non-convex optimal control problems,''
  {\em arXiv preprint arXiv:1804.06539}, 2018.

\bibitem{mao2019convexification}
Y.~Mao, D.~Dueri, M.~Szmuk, and B.~A{\c{c}}{\i}kme{\c{s}}e, ``Convexification
  and real-time optimization for mpc with aerospace applications,'' in {\em
  Handbook of Model Predictive Control}, pp.~335--358, Springer, 2019.

\bibitem{domahidi2013ecos}
A.~Domahidi, E.~Chu, and S.~Boyd, ``Ecos: An socp solver for embedded
  systems,'' in {\em 2013 European Control Conference (ECC)}, pp.~3071--3076,
  IEEE, 2013.

\bibitem{scs}
B.~O'Donoghue, E.~Chu, N.~Parikh, and S.~Boyd, ``{SCS}: Splitting conic solver,
  version 2.1.1.'' \url{https://github.com/cvxgrp/scs}, Nov. 2017.

\bibitem{mao2016successive}
Y.~Mao, M.~Szmuk, and B.~A{\c{c}}{\i}kme{\c{s}}e, ``Successive convexification
  of non-convex optimal control problems and its convergence properties,'' in
  {\em 2016 IEEE 55th Conference on Decision and Control (CDC)},
  pp.~3636--3641, IEEE, 2016.

\bibitem{szmuksuccessive}
M.~Szmuk, {\em Successive Convexification \& High Performance Feedback Control
  for Agile Flight}.
\newblock PhD thesis, University of Washington.

\bibitem{reynolds2019dual}
T.~P. Reynolds, M.~Szmuk, D.~Malyuta, M.~Mesbahi, B.~Acikmese, and J.~M.
  Carson~III, ``Dual quaternion based powered descent guidance with
  state-triggered constraints,'' {\em arXiv preprint arXiv:1904.09248}, 2019.

\bibitem{mao2017successive}
Y.~Mao, D.~Dueri, M.~Szmuk, and B.~A{\c{c}}{\i}kme{\c{s}}e, ``Successive
  convexification of non-convex optimal control problems with state
  constraints,'' {\em IFAC-PapersOnLine}, vol.~50, no.~1, pp.~4063--4069, 2017.

\bibitem{carvalho2013predictive}
A.~Carvalho, Y.~Gao, A.~Gray, H.~E. Tseng, and F.~Borrelli, ``Predictive
  control of an autonomous ground vehicle using an iterative linearization
  approach,'' in {\em 16th International IEEE Conference on Intelligent
  Transportation Systems (ITSC 2013)}, pp.~2335--2340, IEEE, 2013.

\bibitem{gray2013integrated}
A.~J. Gray, {\em An Integrated Framework for Planning and Control of
  Semi-Autonomous Vehicles}.
\newblock PhD thesis, UC Berkeley, 2013.

\bibitem{alrifaee2018real}
B.~Alrifaee and J.~Maczijewski, ``Real-time trajectory optimization for
  autonomous vehicle racing using sequential linearization,'' in {\em 2018 IEEE
  Intelligent Vehicles Symposium (IV)}, pp.~476--483, IEEE, 2018.

\bibitem{mao2018tutorial}
Y.~Mao, M.~Szmuk, and B.~A{\c{c}}{\i}kme{\c{s}}e, ``A tutorial on real-time
  convex optimization based guidance and control for aerospace applications,''
  in {\em 2018 Annual American Control Conference (ACC)}, pp.~2410--2416, IEEE,
  2018.

\bibitem{szmuk2018successive}
M.~Szmuk and B.~Acikmese, ``Successive convexification for 6-dof mars rocket
  powered landing with free-final-time,'' in {\em 2018 AIAA Guidance,
  Navigation, and Control Conference}, p.~0617, 2018.

\bibitem{reynolds2019state}
T.~Reynolds, M.~Szmuk, D.~Malyuta, M.~Mesbahi, B.~Acikmese, and J.~M. Carson,
  ``A state-triggered line of sight constraint for 6-dof powered descent
  guidance problems,'' in {\em AIAA Scitech 2019 Forum}, p.~0924, 2019.

\bibitem{malyuta2019fast}
D.~Malyuta, T.~P. Reynolds, M.~Szmuk, B.~Acikmese, and M.~Mesbahi, ``Fast
  trajectory optimization via successive convexification for spacecraft
  rendezvous with integer constraints,'' {\em arXiv preprint arXiv:1906.04857},
  2019.

\bibitem{ridderhof2019minimum}
J.~Ridderhof and P.~Tsiotras, ``Minimum-fuel powered descent in the presence of
  random disturbances,'' in {\em AIAA Scitech 2019 Forum}, p.~0646, 2019.

\bibitem{vinod2018stochastic}
A.~P. Vinod, S.~Rice, Y.~Mao, M.~M. Oishi, and B.~A{\c{c}}{\i}kme{\c{s}}e,
  ``Stochastic motion planning using successive convexification and
  probabilistic occupancy functions,'' in {\em 2018 IEEE Conference on Decision
  and Control (CDC)}, pp.~4425--4432, IEEE, 2018.

\bibitem{szmuk2017convexification}
M.~Szmuk, C.~A. Pascucci, D.~Dueri, and B.~A{\c{c}}ikme{\c{s}}e,
  ``Convexification and real-time on-board optimization for agile quad-rotor
  maneuvering and obstacle avoidance,'' in {\em 2017 IEEE/RSJ International
  Conference on Intelligent Robots and Systems (IROS)}, pp.~4862--4868, IEEE,
  2017.

\bibitem{augugliaro2012generation}
F.~Augugliaro, A.~P. Schoellig, and R.~D'Andrea, ``Generation of collision-free
  trajectories for a quadrocopter fleet: A sequential convex programming
  approach,'' in {\em 2012 IEEE/RSJ international conference on Intelligent
  Robots and Systems}, pp.~1917--1922, IEEE, 2012.

\bibitem{carvalho2016predictive}
A.~M. Carvalho, {\em Predictive Control under Uncertainty for Safe Autonomous
  Driving: Integrating Data-Driven Forecasts with Control Design}.
\newblock PhD thesis, UC Berkeley, 2016.

\bibitem{zhang2018toward}
Y.~Zhang, H.~Chen, S.~Waslander, T.~Yang, S.~Zhang, G.~Xiong, and K.~Liu,
  ``Toward a more complete, flexible, and safer speed planning for autonomous
  driving via convex optimization,'' {\em Sensors}, vol.~18, no.~7, p.~2185,
  2018.

\bibitem{laumond1998robot}
J.-P. Laumond {\em et~al.}, {\em Robot motion planning and control}, vol.~229.
\newblock Springer, 1998.

\bibitem{kong2015kinematic}
J.~Kong, M.~Pfeiffer, G.~Schildbach, and F.~Borrelli, ``Kinematic and dynamic
  vehicle models for autonomous driving control design,'' in {\em 2015 IEEE
  Intelligent Vehicles Symposium (IV)}, pp.~1094--1099, IEEE, 2015.

\bibitem{rajamani2011vehicle}
R.~Rajamani, {\em Vehicle dynamics and control}.
\newblock Springer Science \& Business Media, 2011.

\bibitem{qian2016model}
X.~Qian, {\em Model predictive control for autonomous and cooperative driving}.
\newblock PhD thesis, 2016.

\bibitem{malyuta2019discretization}
D.~Malyuta, T.~Reynolds, M.~Szmuk, M.~Mesbahi, B.~Acikmese, and J.~M. Carson,
  ``Discretization performance and accuracy analysis for the rocket powered
  descent guidance problem,'' in {\em AIAA Scitech 2019 Forum}, p.~0925, 2019.

\bibitem{antsaklis2006linear}
P.~J. Antsaklis and A.~N. Michel, {\em Linear systems}.
\newblock Springer Science \& Business Media, 2006.

\bibitem{hespanha2018linear}
J.~P. Hespanha, {\em Linear systems theory}.
\newblock Princeton university press, 2018.

\bibitem{gill1981practical}
P.~E. Gill, W.~Murray, and M.~H. Wright, {\em Practical optimization}.
\newblock London: Academic Press, 1981.

\bibitem{ross2018scaling}
I.~Ross, Q.~Gong, M.~Karpenko, and R.~Proulx, ``Scaling and balancing for
  high-performance computation of optimal controls,'' {\em Journal of Guidance,
  Control, and Dynamics}, vol.~41, no.~10, pp.~2086--2097, 2018.

\bibitem{benedikter2019convex}
B.~Benedikter, A.~Zavoli, and G.~Colasurdo, ``A convex approach to rocket
  ascent trajectory optimization,'' in {\em 8th European Conference for
  Aeronautics and Space Sciences (EUCASS)}, 2019.

\bibitem{szmuk2018successive_trig}
M.~Szmuk, T.~P. Reynolds, and B.~Acikmese, ``Successive convexification for
  real-time 6-dof powered descent guidance with state-triggered constraints,''
  {\em arXiv preprint arXiv:1811.10803}, 2018.

\bibitem{cottle1992linear}
R.~W. Cottle, J.-S. Pang, and R.~E. Stone, {\em The linear complementarity
  problem}, vol.~60.
\newblock Siam, 1992.

\bibitem{CVXeliminating}
``Advanced topics, eliminating quadratic forms.''
  \url{http://cvxr.com/cvx/doc/advanced.html}.
\newblock Accessed: 2019-10-18.

\bibitem{dueri2014automated}
D.~Dueri, J.~Zhang, and B.~A{\c{c}}ikmese, ``Automated custom code generation
  for embedded, real-time second order cone programming,'' {\em IFAC
  Proceedings Volumes}, vol.~47, no.~3, pp.~1605--1612, 2014.

\bibitem{kato2018autoware}
S.~Kato, S.~Tokunaga, Y.~Maruyama, S.~Maeda, M.~Hirabayashi, Y.~Kitsukawa,
  A.~Monrroy, T.~Ando, Y.~Fujii, and T.~Azumi, ``Autoware on board: Enabling
  autonomous vehicles with embedded systems,'' in {\em 2018 ACM/IEEE 9th
  International Conference on Cyber-Physical Systems (ICCPS)}, pp.~287--296,
  IEEE, 2018.

\bibitem{aw:autoware_repo}
``Autoware open-source software for self-driving vehicles,'' 2020.

\bibitem{steven}
``Scvxm scpp implementations github repositories.''
  \url{https://github.com/EmbersArc}.
\newblock Accessed: 2019-10-18.

\end{thebibliography}
\end{document}